\documentclass[submission]{FPSAC2027}
\newtheorem{thm}{Theorem}
\newtheorem{lem}{Lemma}
\newtheorem{definition}{Definition}
\newtheorem{expl}{Example}

\DeclareMathOperator{\convspan}{cspan}

\usepackage{mathtools}

\title{Linear algebra over bands}

\author[Yifan Guo and Victoria Schleis]{Yifan Guo\thanks{\href{mailto:yifan_guo3@brown.edu}{yifan\_guo3@brown.edu} Yifan Guo was supported by NSF DMS-2401282.}\addressmark{1} \and Victoria Schleis \thanks{\href{mailto:victoria.m.schleis@durham.ac.uk}{victoria.m.schleis@durham.ac.uk} Victoria Schleis was partially funded by UKRI FLF MR/Y003888/1.}\addressmark{2}}

\address{\addressmark{1} Department of Mathematics, Brown University, Providence, Rhode Island, United States \\ \addressmark{2}Department of Mathematical Sciences, Durham University, Durham, United Kingdom}

\received{\today}

\abstract{We introduce Eckmann-Hilton operators on bands, which allow for a suite of matroidal constructions, in particular, matroid unions and matroid stable intersections. 
This enables the construction of general linear monoids over bands. Geometrically, we show that matroid stable intersection extends the transverse intersection of linear spaces over fields and the stable intersection of tropical linear spaces, lending new insight into intersection theory on bands.}


\usepackage[backend=bibtex]{biblatex}
\begin{document}

\maketitle

\section{Introduction}

Methods from linear algebra are universally useful in all areas of mathematics. Bands, introduced by Baker-Jin-Lorscheid \cite{BakerJinLorscheid}, are a generalization of rings that can be used to describe many settings in algebraic combinatorics, in particular, oriented matroids and tropical linear spaces as in \cite{BakerLorscheid}. Their work builds on earlier studies of matroids with coefficients, developed by Baker-Bowler in \cite{BakerBowler} and originating in work of Dress-Wenzel \cite{dresswenzel}. Beyond matroid theory, bands provide a unified framework to study algebraic geometry in these settings. For this, methods in linear algebra are useful to define and study universal notions of tangent spaces and vector bundles, as well as many more.

Over a band $B$, the analog of a linear space is a $B$-matroid, introduced in \cite{BakerLorscheid}. This leads to challenges in constructing linear maps, 
as they need to be simultaneously general enough to be useful, and defined matroid theoretically. In this extended abstract, we develop a matroidal analog of matrices, allowing us to model linear maps.

The first obstruction is that formal addition over a band $B$ does not yield an element of $B$. 
We overcome this by introducing an \emph{Eckmann-Hilton operator} $\oplus$ in Section \ref{sec: matroids}, which should be thought of as ``generic evaluation" of the formal sum over $B$, see Def.~\ref{def: EH operator}. Past its existence, most of our constructions require an additional compatibility assumption of the formal sum and the Eckmann-Hilton operator, which we call \emph{tolerant}, see Def.~\ref{def: tolerant}. In the following, we require every band $B$ to have an Eckmann-Hilton operator $\oplus$ and to be tolerant. Moreover, we denote the additive inverse of $1$ in $B$ by $\epsilon$.

\begin{thm}\label{mainthm: matroid structure}
Let $B$ be a tolerant band with an Eckmann-Hilton operator $\oplus$, and let $\mu$ and $\nu$ be two $B$-matroids of rank $r$ and $s$ over a ground set $E$, such that $r+s \leq |E|$. Then, the \emph{matroid union} $\mu \cup\nu$, defined coordinate-wise on $Z\in \binom{E}{r+s}$ as $$ \mu\cup \nu (Z) = \bigoplus_{X\sqcup Y = Z} \epsilon^{|\{ y<x \ | \ x\in X,\ y\in Y\}|} \mu(X)\nu(Y)$$ is either identically zero or a $B$-matroid of rank $r+s$ on $E.$ If $\oplus$ is not commutative, the sum is taken in lexicographical order (in the sense of Example \ref{expl: matroid union}). 
\end{thm}
To prove this result, we define several additional matroidal constructions. We show that \emph{direct sums} (see Def.~\ref{def: stable sum}), \emph{reorderings} (see Def.~\ref{def: reordering}), and \emph{contractions along equivalent elements} (see Def.~\ref{def: contraction}) of $B$-matroids still satisfy the Grassmann-Pl\"ucker relations.

Using these matroidal methods allows us to give a well-defined notion of linear algebra over bands in Section \ref{sec: matrices and bimatroids}. 
As a stand-in for $n\times n$ matrices, we consider a $B$-matroid of rank $n$ on $[2n]$ with a fixed basis $[n]$. We call such a matroid a \textit{$B$-bimatroid}. Unlike $B$-matroids or $B$-matrices, these have a well-defined and well-behaved multiplication. We call the  set of full-rank $n\times n$ $B$-bimatroids the \emph{general linear monoid} $GL_B(n)$.
\begin{thm}\label{mainthm: bimatroid product}
Let $B$ be a tolerant band with an Eckmann-Hilton operator, and let $\mathsf{M}$ and $\mathsf{N}$ be full-rank $n\times n$-$B$-bimatroids, then the product $\mathsf{M} \cdot \mathsf{N}$ is a full-rank $n \times n$-$B$-bimatroid. Further, the bimatroid product is associative; and the unit bimatroid induced by the identity matrix is a unit of this product. In other words, $GL_B(n)$ is a monoid.
\end{thm}

In Section \ref{sec: intersections},  we shift to a geometric perspective. In tropical geometry, many statements on valuated matroids can be rewritten as statements about tropical linear spaces. These reformulations often allow for the proof of further statements on bands using geometric intuition. 

We first show that the $B$-matroid stable intersection $\mu \cap \nu$ (see Def.~\ref{def: matroid stable intersection}) is 
dual to the matroid union of the dual matroids $\mu^\perp \cup \nu^\perp$. Using the Eckmann-Hilton operator and Baker-Bowler's \emph{cocircuits} $\mathcal{C}^{\ast}$ \cite{BakerBowler}, we define the \emph{$B$-linear space} $\mathcal{W}^{\ast}$ (see Def. \ref{def: convex covectors}), which is cryptomorphic to a $B$-matroid for any band $B$ in which all nonzero elements are units. Such bands are called \textit{idylls} in \cite{BakerJinLorscheid}. 

\begin{thm} \label{mainthm: geometry of B-linear spaces}
Let $B$ be a tolerant idyll with an Eckmann-Hilton operator. Let $\mu$ and $\nu$ be $B$-matroids with nonzero matroid stable intersection $\mu \cap \nu$. Then, the $B$-linear space of $\mu \cap \nu$ is contained in the set-theoretic intersection of $B$-linear spaces for $\mu$ and $\nu$, i.e.,
    \[\mathcal{W}^{\ast\mu \cap \nu} \subset \mathcal{W}^{\ast \mu}\cap \mathcal{W}^{\ast \nu},\]
with equality under additional assumptions on the transversality of the intersection, convexity of the linear spaces, and fusion property of the band $B$.
\end{thm}

In particular, for tropical geometers,  Theorems \ref{mainthm: matroid structure} and \ref{mainthm: geometry of B-linear spaces} provide alternative proofs 
to results on stable sums and stable intersections of valuated matroids 
in \cite{giansiracusa2017grassmannalgebramatroids,speyer2004tropicallinearspaces}. For oriented matroid theorists, we fill in a missing step in the existing proof of oriented matroid union in \cite{LAWRENCE1981183} in  Theorem \ref{mainthm: matroid structure}, and introduce new concepts in Theorem \ref{mainthm: bimatroid product}. 

\section{Matroid unions over bands}\label{sec: matroids}
In this section, we study matroids over bands with Eckmann-Hilton operators. We begin by recalling the theory of bands and introduce the Eckmann-Hilton operator, an additional formal structure on bands that can be interpreted as a \emph{generic evaluation map}. We then use this structure to generalize matroid union. 

\subsection{Tolerant bands with Eckmann-Hilton operators}
In this section we recall bands, which were introduced in \cite{BakerJinLorscheid}, and then proceed to define the Eckmann-Hilton operator and what it means for bands to be tolerant.

\begin{definition}
    A \emph{band} $B$ is a commutative 
    monoid  $(B,\cdot,1)$, together with a set $N_B$,  called the \emph{null set},  consisting of formal sums $\sum a_i$ of elements $a_i \in B$, such that 
    \begin{description}
        \item[(Zero element)] there exists an element $0\in B$ that satisfies $a\cdot 0 = 0$ for all $a\in B$;
        \item[(Null ideal)] we have $0 \in N_B$,  $a \cdot x \in N_B$, and $x+y \in N_B$ for all $a \in B$ and $x,y \in N_B$; and
        \item[(Additive inverse)] for every $a\in B$, there exists a unique $b\in B$ with $a+b \in N_B$.
    \end{description}In every band, there exists a distinguished element $\epsilon$ such that $1 +\epsilon \in N_B$. 
\end{definition}

Bands form a well-defined and reasonable category,
as elaborated on in \cite{BakerJinLorscheid}. Furthermore, many structures of interest that arise in combinatorial algebraic geometry can be described as structures over bands. 

\begin{expl}\label{ex:bands}
Every ring $R$ is a band, where $N_R = \{ \sum a_i \ | \ \sum a_i = 0 \ \text{  evaluated in } R\}$. 
    The \emph{Krasner band} $\mathbb{K} = \{0,\, 1\}$ is defined with $N_{\mathbb{K}} = \{\sum_{i=1}^{k}1 \ | \ k\geq 2\} \cup \{0\}.$ The \emph{tropical band} $\mathbb{T}=\mathbb{R}_{\geq 0}$ is defined with the null set consisting of $0$ and all formal sums that attain the maximum of their summands at least twice. The \emph{sign band} $\mathbb{S} = \{0, 1, \epsilon\}$ is defined with the null set consisting of $0$ and formal sums that contain at least one $1$ and one $\epsilon$.
\end{expl}

When studying linear algebra over bands, an immediate problem is that formal sums of elements of $B$ are not elements of $B$. This is problematic when considering matrix multiplication. To this end, we consider bands with an additional structure:

\begin{definition}\label{def: EH operator}
Let $B$ be a band. An \emph{Eckmann-Hilton (EH) operator} is an associative binary operation $\oplus: B\times B \rightarrow B$
that satisfies 
\begin{description}
    \item[(EH1)] that the zero-element of $B$ is the additive identity, $0 \oplus x = x \oplus 0 = x$ for all $x\in B$;
    \item[(EH2)] double distributivity, that is, $(a \oplus b) (c \oplus d) = ac \oplus bc \oplus ad \oplus bd$ for all $a,b,c,d \in B$;
    \item[(EH3)] for all $a, b \in B^m $ that  if $\sum_{i =1}^m a_i \in N_B$ and   $\sum_{i =1}^m b_i \in N_B$, then $\sum_{i = 1}^m(a_i \oplus b_i) \in N_B$.
\end{description}
\end{definition}

An EH operator should be thought of as ``choosing a generic representative" to evaluate the formal sum.  The axiom (EH3) should be thought of as allowing $\oplus$ to ``commute'' with the formal sum. In particular, this partially satisfies the hypothesis of the topological Eckmann–Hilton argument, which this definition is named after.

\begin{expl}\label{ex: EH operators}
    Most bands of interest have EH operators.  
    \begin{itemize}
    \item For a field $K$, the EH operator is the evaluation of the sum using standard addition of $K$. 
  \item The EH operator over the Krasner band $\mathbb{K}$ is induced by $1 \oplus 1 = 1$. 
            \item The EH operator over the tropical band $\mathbb{T}$ is the standard tropical addition, given by $a \oplus b = \max(a,b)$ for all $a, b \in \mathbb{T}$. 
            This can be generalized for a broader class of
            bands that includes the tropical real hyperfield, tropical complex hyperfield, and tropical phase hyperfield in the sense of \cite{anderson2018vectorsmatroidstracts}.
            \item Over the band of signs $\mathbb{S}$ there are two EH operators, which we denote by $\oplus_{>}$ and $\oplus_{<}$. They are induced by $1 \oplus_{>} \epsilon = 1$, and $ \epsilon\oplus_{>} 1 = \epsilon$, and respectively by $1 \oplus_{<} \epsilon = \epsilon$, and $\epsilon \oplus_{<} 1 = 1$. Note that this is an example of a non-commutative EH operator. This can also be generalized for all \emph{fusion idylls with the inflation property}, as defined in \cite[Def.~1.9]{BakerJinLorscheid} and  \cite[Def.~6.9]{anderson2018vectorsmatroidstracts}.
          
    \end{itemize}
\end{expl}

Bands allow for a unified study of matroids, defined by \emph{Grassmann-Pl\"ucker relations}.

\begin{definition}[see {\cite{BakerJinLorscheid}}]\label{def: bmatroid}
    A \emph{$B$-matroid of rank $r$ on a totally ordered ground set $E$} is an equivalence class of functions  $\mu: \tbinom{E}{r} \rightarrow B$, such that $\mu(I)\in B^\times$ for some $I\in \tbinom{E}{r}$, and the following \emph{Grassmann-Pl\"ucker relations} are satisfied for any $I\in \tbinom{E}{r+1}$ and  $J\in \tbinom{E}{r-1}$,
    \begin{equation}\label{pluecker}
     \sum_{k=0}^{r} \epsilon^{|\{e \in I \cup J \ \mid\ e < i_k\}|} \mu(I\setminus \{i_k\})\mu(J \cup \{i_k\}) \in N_B.
    \end{equation}
    We consider $\mu$ up to projective equivalence, i.e., with $\mu\sim a\cdot \mu \ $ for all $a \in B^\times$.
\end{definition}

Many familiar classes of matroids can be recovered as $B$-matroids. For a field $K$, a $K$-matroid of rank $r$ on $[n]$ is an  $r$-dimensional subspace of $K^n$. Classical matroids are precisely matroids over $\mathbb{K}$, oriented matroids are $\mathbb{S}$-matroids, and valuated matroids are $\mathbb{T}$-matroids.
\begin{definition}\label{def: tolerant}
    A band $B$ with EH operator $\oplus$ is  \emph{tolerant} if any $a_i, b_i, c_i, d_i, S \in B$ such that
\begin{align*}
    \sum_{i=1}^na_ic_i \in N_B, \ \  \ \ 
    S+ \sum_{i=1}^nb_ic_i \in N_B, \ \  \ \ 
    \epsilon S +  \sum_{i=1}^na_id_i \in N_B, \ \ \text{ and } \ \  
     \sum_{i=1}^nb_id_i \in N_B,
\end{align*}
also satisfy that
$$\sum_{i=1}^n (a_i\oplus b_i)(c_i\oplus d_i) \in N_B.$$
\end{definition}
All  bands we discussed in Example \ref{ex: EH operators} are tolerant with the EH operators given thereafter. 
However, there are also well-known examples which are not tolerant with respect to their canonical EH operator, one being the classical triangular hyperfield:
\begin{expl}
    The triangular hyperfield $\mathbb{T}_2$ where $\oplus$ is addition is not tolerant. For example, 
    \[(1\oplus 1)(1\oplus 1)+(1\oplus 1)(1\oplus 1)+(1\oplus 2)(2\oplus 1) \notin N_B\]
    since there is no triangle with side lengths $4,4,$ and $9$. However, with $S = 2.5$, we have 
    \begin{align*}
        &1\cdot 1 + 1\cdot 1+1\cdot 2 \in N_B \ \ \ \ \ \  
        S+ 1\cdot 1 + 1\cdot 1+2\cdot 2 \in N_B\\
        &S+ 1\cdot 1 + 1\cdot 1+1\cdot 1 \in N_B
        \ \ \ \ \ \  1\cdot 1 + 1\cdot 1+2\cdot 1 \in N_B.
    \end{align*}
\end{expl}

\subsection{Matroid unions over tolerant bands with EH operators}
The main goal of this section is to give a sketch of the proof of Theorem \ref{mainthm: matroid structure}, which uses an inductive technique inspired by \cite{LAWRENCE1981183}.
The full proof is long and technical, and exceeds the scope of this extended abstract. Instead, we outline the main steps, which involve showing that several matroidal constructions are again $B$-matroids.

\begin{definition}\label{def: stable sum} 
    Let $B$ be a band. Let $\mu$ be a matroid of rank $r$ on a ground set $E$, and let $\nu$ be a matroid of rank $s$ on a ground set $F$. Then, their \emph{direct sum} $\mu \widehat\oplus \nu$ is the map
    \begin{align*}
    \mu \widehat\oplus \nu : \tbinom{E \sqcup F}{r+s} & \rightarrow B; \ \  \ Z \mapsto \mu(Z\cap E) \nu(Z\cap F).
    \end{align*}
    For the ground set $E\sqcup F$ of $\mu\widehat\oplus \nu$, we fix the order $e< f$ for any $e\in E$ and $f\in F$, and inherit the order on $E$ and $F$ otherwise. 
\end{definition}
By a direct computation, it can be shown that the direct sum $\mu \widehat\oplus \nu$  is a $B$-matroid. Note that the direct sum with a uniform matroid is a natural generalization of elongation of valuated matroids as defined in \cite[Sec.~5.3]{giansiracusa2017grassmannalgebramatroids}.

Since the EH operator is not necessarily commutative, it is sometimes necessary to reorder elements of the ground set, flipping signs.

\begin{definition}\label{def: reordering}
    Let $B$ be a band. For a $B$-matroid $\mu :\tbinom{E}{r} \rightarrow B$  and two consecutive elements $a<b \in E$, we define a \emph{reordering} of $\mu$ to be the new function 
    \begin{align*}
        \mu_{b<a} :\tbinom{E'}{r} &\rightarrow B,
        \ \ \ \   \mu_{b<a} (Z)  = \begin{cases} 
      \epsilon \mu(Z) & \text{ if } a,b \in Z \\
      \mu(Z) & \text{otherwise, }  
   \end{cases}
    \end{align*}
    where $E' = E$ as a set, with identical ordering, except for the order of $a$ and $b$ switched.
\end{definition}
As before, it can be shown by a computation that $\mu_{b<a}$ is a $B$-matroid.

\begin{definition}\label{def: contraction}
    Let $\mu: \tbinom{E}{r}\rightarrow B$ be a $B$-matroid, and let $u < v \in E$ be two distinct ground set elements, adjacent in the total order of $E$. Consider $E' \coloneqq E/(u \sim v)$, and set $w \coloneqq [u] = [v]$ in $E'$.  
    The order on $E'$ is inherited from the order on $E$ by setting $x <_{E'} w$ if and only if $x<_{E} u$.
    
    We define the \emph{contraction of $\mu$ identifying $u$ and $v$} as
    \begin{equation}
        \xi : \tbinom{E'}{r} \rightarrow B, \ \  \ \  \ \  
         \xi(Z) = \begin{cases}
            \mu(Z) & \text{if } w \notin Z \\
            \mu((Z\setminus w) \cup u) \oplus \mu((Z\setminus w) \cup v) & \text{if } w\in Z.
        \end{cases}
    \end{equation}
\end{definition}

\begin{thm}\label{prop: contraction} Suppose $B$ is tolerant band with EH operator $\oplus$.
     Let $\mu: \tbinom{E}{r}\rightarrow B$ be a $B$-matroid, let $u\neq v \in E$ be adjacent terms in the total order of $E$, and let $\xi$ be the contraction of $\mu$ identifying $u$ and $v$ as defined above. Then $\xi$ satisfies the  Grassmann-Pl\"ucker relations.
\end{thm}
\begin{proof}[Proof sketch]
    The proof is lengthy and proceeds along a careful case-by case analysis, distinguished by whether the exchanged element is contained in a set with $w$. The most challenging case is when $w$ is contained in both $I$ and $J$ in the Grassmann-Pl\"ucker relation of $\xi$ given in \eqref{pluecker}. In this case, we utilizes the four Pl\"ucker relations of $\mu$ involving $I\setminus w \cup u, I\setminus w \cup v$ with $J\setminus w \cup u, J\setminus w \cup v$. Unlike the direct sum in Definition~\ref{def: stable sum}, this case utilizes the assumption that the band is tolerant.
\end{proof}
It is possible that the contraction of a $B$-matroid $\mu$ identifying $u$ and $v$ is identically zero. In this case, there exists a $B$-matroid $\xi$ of smaller rank that can take its place. It is defined as
     \begin{equation}\label{eq: lower rank contract}
         \xi: \tbinom{E'}{r-1} \rightarrow B,  \ \ \   \xi(U) = \begin{cases}
            \mu(U \cup v ) & \text{if } w \notin U \\
            \mu((U\backslash w )\cup\{u,v\}) & \text{if } w\in U.
        \end{cases}
     \end{equation}
\color{black}

We are now ready to give a proof of Theorem \ref{mainthm: matroid structure}, showing that the matroid union $\mu \cup \nu$ satisfies the Pl\"ucker relations, that is, is either identically zero or a $B$-matroid.

\begin{proof}[Proof  of Theorem~\ref{mainthm: matroid structure}]
    Consider the direct sum $\mu \widehat\oplus \nu$ on the disjoint union of ground sets $E\sqcup E$. It can be shown that this function satisfies the Pl\"ucker relations. To distinguish between the two copies of $E$, denote elements of the first copy by $1, \ldots,l$ and elements of second copy by $1', \ldots,l'$. Consider a reordering $\mu \widehat\oplus \nu_{lex}$ of $\mu \widehat\oplus \nu$ arising by sequentially reordering $E\sqcup E$ from $1< \ldots<l<1'< \ldots<l'$ to $1<1'<2<2'<\cdots< l< l'$. This is still a $B$-matroid by the discussion after Definition \ref{def: reordering}. 
    
    Starting with $\mu \widehat\oplus \nu_{lex}$, we sequentially apply matroid contraction by identifying $i$ with $i'$ for $i=l,l-1, \ldots, 1$ in decreasing order in $E$. Denote the resulting function after contracting up to $j$ by $\xi_j$. All $\xi_j$ still satisfy the the Pl\"ucker relations by Theorem \ref{prop: contraction}.
    
    We verify that $\xi_1$ is equal to $\mu \cup \nu$: 
    For any $Z = \{z_1 < \cdots < z_{r+s}\} \in \tbinom{E}{r+s}$, let $u_i, v_i$ be the two copies of $z_i$ in $E\sqcup E$ before the contractions. We write elements of the set in strictly ascending order.  Inductively writing $\xi_{z_i}$ in terms of $\xi_{z_{i+1}}$, we obtain
    \begin{align*}
        \xi_1(z_1,\cdots,z_{r+s}) 
        \ \ = \ \  &\xi_{z_2}(u_1,\cdots ,z_{r+s}) \oplus \xi_{z_2}(v_1,\cdots, z_{r+s}) \\
        \overset{(Ind.)}{=}    &\mu \widehat\oplus \nu_{lex}(u_1,u_2,\cdots, u_{r+s}) \oplus \cdots \oplus \mu \widehat\oplus \nu_{lex}(v_1,v_2 ,\cdots , v_{r+s})\\
        = \ \  &\bigoplus_{X\sqcup Y = Z \text{ lex.}} \epsilon^{|\{y \leq x \mid x \in X, y \in Y\}|}\mu(X)\nu(Y)\\
        = \ \ & \mu \cup \nu (Z).
    \end{align*}

    Thus, $\mu \cup \nu$ satisfies the Pl\"ucker relations. This means that the matroid union is either identically zero or a $B$-matroid itself.
\end{proof}

    In the case that $\mu \cup \nu$ is identically zero, there may be a $B$-matroid of a smaller rank that can be constructed by using the lower rank contraction in Equation \eqref{eq: lower rank contract},
    whenever contracting by identifying $i\sim i'$ results in a zero function. When $B$ is $\mathbb{K}$ or $\mathbb{S}$, this $B$-matroid of a potentially smaller rank respectively coincides with the usual notion of (oriented) matroid union.

    If the EH operator $\oplus$ is commutative, i.e., the semiring arising from the EH operator is idempotent, then $\mu\cup \nu = \epsilon^{rs} \nu\cup \mu \sim \nu\cup \mu$ since 
    $$|\{ y<x \ | \ x\in X,\ y\in Y\}| = |X||Y| - |\{ x<y \ | \ x\in X,\ y\in Y\}|.$$

\begin{expl}\label{expl: matroid union}
Let $\mu$ be a $\mathbb{T}$-matroid of rank 1 on $[4]$ where $\mu(1) = \mu(3)=1, \mu(2) = 2, \mu(4) = 0,$
    and let $\nu$ be a $\mathbb{T}$-matroid of rank 2 on $[4]$ given by 
    \[\nu(23)=\nu(24)=\nu(34)=1, \nu(12)=\nu(13)=\nu(14)=0.\]
Then, the matroid union $\mu \cup \nu$ is the rank 3 matroid on $[4]$ with
    \[\mu \cup \nu (234) = \mu(2)\nu(34) \oplus \mu(3)\nu(24) \oplus \mu(4)\nu(23)= \max\{2, 1, 0\} = 2,\]
    and similarly $\mu \cup \nu(134)=\mu \cup \nu(124) = \mu \cup \nu(123) = 1$.
\end{expl}

\section{Bimatroid products over bands}\label{sec: matrices and bimatroids}
For a band $B$, we now construct the matroidal analog of matrices. Although the $B$-matroid defined by taking the maximal minors of a matrix with entries in $B$ using the EH operator can yield a first approximation, the data contained therein is not sufficient for many constructions. Thus, we define a larger $B$-matroid that captures the information of all minors of the matrix, rather than just the maximal ones.
\begin{definition}
    A \emph{bimatroid} $\mathsf{M}$ over a band $B$ on rows $E$ and columns $F$ is a $B$-matroid $\mathsf{M}: \binom{E\sqcup F}{|E|}\rightarrow B$ where $\mathsf{M}(E)  = 1$. In shorthand, we say that $\mathsf{M}$ is an $|E|\times |F|$-bimatroid.
    
    Equivalently, we  express $\mathsf{M}$ as a bivalued map $\mathsf{M}(R,C) \coloneqq \mathsf{M}((E-R )\cup C),$ defined for all $R \subset E, C \subset F$ with $|R|=|C|$.
\end{definition}
$\mathbb{K}$-bimatroids are precisely the \emph{bimatroids} introduced by Kung and Schrijver  \cite{Kung, Schrijver}, i.e. matchings on bipartite graphs. Over $\mathbb{T}$, we recover \emph{valuated bimatroids}, introduced by Murota \cite{MUROTA}.

A natural example, and the motivation of the definition of bimatroids, is that minors of $B$-matrices induce $B$-bimatroids. Assume $m\geq n$ and consider a matrix $A \in B^{n \times m}$. Appending an $n\times n$ identity matrix $(I|A)$ yields a  $n \times (n+m)$ matrix. Taking the maximal minors (with respect to the EH operator) of this matrix yields a bimatroid. The fact that this is a matroid is not straightforward, and can be  proven using matroid union.

\begin{definition}\label{def: product bimatroids}
Let  $\mathsf{M}$ be an $n\times m$-bimatroid and $\mathsf{N}$ be an $m\times l$-bimatroid over a band $B$. Their \emph{product} $\mathsf{M}\cdot \mathsf{N}$ is defined by
\[\mathsf{M} \cdot \mathsf{N} (I,K) = \bigoplus_{J\in \tbinom{[m]}{d}} \epsilon^{|\{y<x \mid x \in J, y \in [m]-J\}|}\mathsf{M}(I,J)\mathsf{N}(J,K)\]
for $(I, K)\in\tbinom{[n]}{d}\times\tbinom{[l]}{d}$ and any $0 \leq d \leq \min(n,l)$, where $J$ is summed in lexicographical order.
\end{definition}

When $\mathsf{M},\mathsf{N}$ are given by all minors of matrices $A,A'$ over a field, the bimatroid product $\mathsf{M}\cdot \mathsf{N}$ coincides with all minors of the matrix product $A A'$, due to the Cauchy–Binet formula. Over $\mathbb{T}$, Definition \ref{def: product bimatroids} recovers ``composition"  as defined in \cite[Pg.~107]{Frenk}.

\begin{expl}\label{expl: bimatroid product}
Let $\mathsf{M}$ be a $\mathbb{T}$-bimatroid on row $\{1\}$ and columns $\{2,3\}$ given by
    \[\mathsf{M}(\emptyset, \emptyset) = \mathsf{M}(1,3) = 1, \mathsf{M}(1,2) = 2,\]
    and let $\mathsf{N}$ be a $\mathbb{T}$-bimatroid of on rows $\{2,3\}$ and column $\{4\}$ given by 
    \[\mathsf{N}(\emptyset,\emptyset) = \mathsf{N}(2,4) = \mathsf{N}(3,4) = 1.\]
Then, the bimatroid product $\mathsf{M} \cdot \mathsf{N}$ is given by
    \[\mathsf{M} \cdot \mathsf{N} (1,4) = \mathsf{M}(1,2) \mathsf{N}(2,4) \oplus \mathsf{M}(1,3) \mathsf{N}(3,4) = \max\{2, 1\} = 2,\]
    and similarly $\mathsf{M} \cdot \mathsf{N}(\emptyset,\emptyset) = 1$.
\end{expl}
To show that the bimatroid product is still a bimatroid over bands, we need two auxiliary matroid constructions.  
\begin{definition}\label{lem: auxiliary}
    Let $\mu: \tbinom{E}{r} \rightarrow B$ be a $B$-matroid, $S$ a finite set, and $T \subset E$ with $|T|<r$.
     
    \begin{enumerate}
        \item Its \emph{loop extension} by $S$ is defined to be $\mu \oplus S(X) := \mu(X)$ if $X \subset E$ and $\mu \oplus S(X):= 0$ otherwise, for any $X \in \tbinom{E \sqcup S}{r}$.
  
   \item Its \emph{generic contraction} by $T$ is defined to be  $\mu/T (X) := \mu(X\cup T)$
    for any $X \in \tbinom{E \setminus T}{r-|T|}$.
    \end{enumerate}
   
\end{definition}

\begin{thm}\label{prop: mult B bimatroid} Let $B$ be a tolerant band with an EH operator. Let $\mathsf{M}$ be an $n\times m$-bimatroid and $\mathsf{N}$ be an $m\times l$-bimatroid over $B$. Their  product $\mathsf{M}\cdot \mathsf{N}$  is a $B$-bimatroid,  that is, $B$-bimatroid multiplication is well-defined.
\end{thm}
\begin{proof}[Proof sketch]
First,  show that both constructions in Definition \ref{lem: auxiliary}  yield a $B$-matroid. Now,
note that 
the bimatroid product $\mathsf{M} \cdot \mathsf{N}$
can be written as
\begin{equation} \label{bimatroid decomp}
    \mathsf{M} \cdot \mathsf{N}= \left((\mathsf{M} \oplus [l]) \cup ([n] \oplus \mathsf{N} ) \right)\, /\, [m].
\end{equation}
Since the result of every one of these constructions satisfy the Pl\"ucker relations, so does the bimatroid product. The condition on the first coordinate follows from the fact that 
$\mathsf{M} \cdot\mathsf{N}(\emptyset, \emptyset) = \mathsf{M} (\emptyset, \emptyset)\mathsf{N}(\emptyset, \emptyset) =  1.$
\end{proof}

We say that an $n\times n$ bimatroid $\mathsf{M}$ is of \emph{full rank} if $\mathsf{M}([n], [n])\in B^\times$, as $\mathsf{M}([n], [n])$ is the analog of the $n\times n$-minor. Then, products of full rank $n\times n$-bimatroids are again full rank, as  $\mathsf{M} \cdot \mathsf{N}([n], [n]) = \mathsf{M}([n], [n])\mathsf{N}([n], [n]) \in B^\times,$ and associativity of the bimatroid product follows from the associativity of the EH operator. Combined, we obtain Theorem \ref{mainthm: bimatroid product}. 

\begin{expl} We demonstrate Equation \eqref{bimatroid decomp} on the  bimatroid product $\mathsf{M} \cdot \mathsf{N}$ from Example \ref{expl: bimatroid product}. We have that
$\mathsf{M} \oplus \{4\} = \mu$ and $\{1\} \oplus \mathsf{N} = \nu$ for $\mu$ and $\nu$ as given in Example \ref{expl: matroid union}. We observe that the map $\left(\mu \cup \nu \right)/\{2,3\}$ is given by $\big((\mu \cup \nu )/\{2,3\}\big)(1) = 1 $ and $ \big((\mu \cup \nu )/\{2,3\}\big)(4) = 2$. This is precisely the matroid form of $\mathsf{M} \cdot \mathsf{N}$.
\end{expl}
\section{Intersections of linear spaces over bands}\label{sec: intersections}

We first formally define \emph{matroid stable intersection} over bands, then relate this construction to both $B$-matroid union and intersections of $B$-linear spaces.

\begin{definition}\label{def: matroid stable intersection}
    Let $B$ be a band with EH operator $\oplus$. For two $B$-matroids $\mu$ of rank $r$ and $\nu$ of rank $s$ on $E$ with $r+s \geq |E|$, we define the \emph{matroid stable intersection} $\mu \cap \nu$ to be the map
    \begin{align*}
    \mu \cap \nu : \tbinom{E}{r+s-|E|} & \rightarrow B;  \ \ \  \ \  \ Z \mapsto \bigoplus_{X\cap Y = Z} \epsilon^{|\{ x<y \ | \ x\in X-Z,\ y\in Y-Z\}|} \mu(X)\nu(Y)
    \end{align*}
    where the sum is taken in reverse lexicographical ordering on $E-Z$. 
\end{definition}

To link stable intersection to matroid union, we recall from \cite{BakerLorscheid} that the \textit{dual} of a $B$-matroid $\mu:\tbinom{E}{r} \rightarrow B$ is the $B$-matroid $\mu^\perp:\tbinom{E}{|E|-r} \rightarrow B$ given by $\mu^\perp(Z) = \mu(E\setminus Z)$.
\begin{expl}\label{ex matroid intersection}
    We now consider the dual $\mathbb{T}$-matroids of the matroids in Example \ref{expl: matroid union}. The $\mathbb{T}$- matroid $\mu^\perp$ has  rank 3 on $[4]$ and is given by
    \[\mu^\perp(234)=\mu^\perp(124)=1, \mu^\perp(134)=2,  \ \ \text{ and } \ \ \mu^\perp(123)=0,\]
    and the $\mathbb{T}$-matroid $\nu^\perp$ has rank 2 on $[4]$, given as
    \[\nu^\perp(14) = \nu^\perp(13)= \nu^\perp(12)=1, \  \  \text{ and } \ \  \nu^\perp(34) =\nu^\perp(24) =\nu^\perp(23) = 0.\]
Their matroid stable intersection $\mu^\perp \cap \nu^\perp$ is given by
    \[\mu^\perp \cap \nu^\perp (1) = \mu(134)\nu(12)\oplus\mu(124)\nu(13)\oplus\mu(123)\nu(14) = \max\{2,1,0\} =2, \]
    and similarly $\mu^\perp \cap \nu^\perp(2)=\mu^\perp \cap \nu^\perp(3) = \mu^\perp \cap \nu^\perp(4) = 1$, which
    is the dual matroid of $\mu \cup \nu$.
    \end{expl}
\begin{thm} \label{thm: matroid intersection}
    Let $\mu^\perp$ and $\nu^\perp$ denote the dual of matroids $\mu$ and $\nu$ over a band $B$ with an EH operator. Then, \[\mu \cap \nu = (\mu^\perp \cup \nu^\perp)^\perp.\]
    Moreover, if $B$ is tolerant, then $\mu \cap \nu$ is either identically zero or a $B$-matroid.
\end{thm}
 The expression follows by algebraic manipulation, and the last statement follows from Theorem \ref{mainthm: matroid structure}. Over $\mathbb{T}$, combining Theorems \ref{mainthm: matroid structure} and \ref{thm: matroid intersection} yields an alternative proof for \cite[Prop.~3.1]{speyer2004tropicallinearspaces}.

We now shift to a geometric perspective. For fields, matroid stable intersection yields the shuffle formula for vector space intersection in \cite{WOHLHART201359}. In tropical geometry, matroid stable intersection can be seen as the stable intersection of tropical linear spaces in \cite{speyer2004tropicallinearspaces}. We generalize this idea to bands. 

When $B$ is an idyll, $B$-matroids can be equivalently defined using the cryptomorphism with cocircuits, which allow us to view $B$-matroids as linear spaces over $B$.

\begin{definition}[{\cite[Def.~4.6]{BakerBowler}}] \label{def: cocircuits covectors}
    For a $B$-matroid $\mu :\tbinom{E}{r} \rightarrow B$ over $E=[n]$, the set of \emph{cocircuits} is $\mathcal{C}^{\ast\mu} = \{\alpha C_S^{\ast\mu}  \ |\  S \in \tbinom{E}{r-1}, \  \alpha \in B^\times \}\setminus \{(0, \ldots, 0\},$
    where 
    $C_S^{\ast\mu} \in B^n$ is given as
    \[C_S^{\ast\mu}(i) = \begin{cases} 
      0 & i \in S \\
      \epsilon^{|\{j > i \ |\  j \in S\}|} \mu(S \cup i) & i \notin S. 
   \end{cases}\]
   \end{definition}
Extending the tropical setting, we define $B$-linear spaces as the span of cocircuits using the Eckmann-Hilton operator.
\begin{definition}\label{def: convex covectors}
We say that $X \in B^n$ is \emph{in the  convex span} of a set of vectors $S = \{Y_k \in B^n \mid k \in [l]\}$, denoted by $\convspan (S)$, if
 there exist $b_1, \ldots, b_l\in B$ such that
     \[X(i) = \bigoplus_{k=1}^l b_kY_k(i) \ \ \  \ \  \text{ for all $i \in [n]$}.\]
For a $B$-matroid $\mu$ of rank $r$ over $[n]$, the \emph{$B$-linear space} is the convex span of its cocircuits, 
    $$\mathcal{W}^\ast:= \convspan(\mathcal{C}^{\ast\mu}).$$
\end{definition}

  We say that a $B$-linear space is \emph{convex} if the $B$-linear space defined as above coincides with the set of covectors as defined in \cite[Def.~2.15]{anderson2018vectorsmatroidstracts}. We say that the band $B$ itself is \emph{convex} if every $B$-linear space is convex. Many key examples of bands are convex --   $\mathbb{T}$ due to  \cite[Thm.~3.8]{MT01}, and both  $\mathbb{K}$ and $\mathbb{S}$ because of how their null sets are defined.

\begin{proof}[Proof sketch of Theorem \ref{mainthm: geometry of B-linear spaces}]
We first expresses the cocircuits of the matroid stable intersection in terms of the cocircuits of the constituent matroids.

\begin{lem}\label{lem: cocircuits of matroid intersection}
    The cocircuit of $\mu\cap \nu$ associated with $S$ can be expressed as
   \[C_S^{\ast\mu\cap \nu}(i) =
      \bigoplus_{X\cap Y = S, X \cup Y = E-i} \epsilon^{|\{x < y \ | \ x\in X-S, \ y \in Y-S\}|+|\{(i< x \ | \ x\in E)\}|} C_X^{\ast\mu}(i)C_Y^{\ast\nu}(i) 
   \]
\end{lem}
To prove that $C_S^{\ast\mu\cap \nu}$ is in the convex span of cocircuits of $\mu$ (or $\nu$ respectively), we show that the coefficient of each $C_X^{\ast\mu}$ in the linear expression is independent of $i$. 
The inclusion $\mathcal{W}^{\ast\mu \cap \nu} \subset \mathcal{W}^{\ast \mu}\cap \mathcal{W}^{\ast \nu}$ then follows from the fact that $B$-linear spaces are closed under the EH operator $\oplus$. Under an additional transversality condition, the reverse inclusion can be shown for convex $B$-linear spaces over fusion band $B$ by manipulating relations on covectors of imposed by circuits of $\mu$ and $\nu$.
\end{proof}

In general, the additional condition needed for $\mathcal{W}^{\ast\mu \cap \nu} = \mathcal{W}^{\ast \mu}\cap \mathcal{W}^{\ast \nu}$ can be seen as requiring the intersection to be ``transversal". Over $\mathbb{K}$, this coincides with \textit{transverse matroids} in \cite{speyer2004tropicallinearspaces}. Over $\mathbb{T}$, it corresponds to \text{transverse intersections} of tropical linear spaces in \cite[Def.~3.4.9]{MaclaganSturmfels2015}, thus providing an alternative proof for \cite[Prop.~3.4]{speyer2004tropicallinearspaces}.

\begin{expl}
  The $\mathbb{T}$-linear spaces of $\mu^\perp$ and $\nu^\perp$ defined in Example \ref{ex matroid intersection} can be described as  
    \begin{align*}
        \mathcal{W}^{\ast \mu^\perp} &= \{v \in B^4 \mid v_1=v_3\geq 2v_2 \text{ or }v_1 = 2v_2\geq v_3 \text{ or } 2v_2 = v_3\geq v_1\}, \text{ and}\\ 
        \mathcal{W}^{\ast \nu^\perp} &= \{v \in B^4 \mid v_2=v_3=v_4\},
    \end{align*}
    with set-theoretic intersection
    \[\mathcal{W}^{\ast \mu^\perp}\cap \mathcal{W}^{\ast \nu^\perp} = \{v \in B^4 \mid v_1=2v_2=2v_3 = 2v_4\},\]
    which indeed equals $\mathcal{W}^{\ast\mu^\perp \cap \nu^\perp}.$
\end{expl}

\acknowledgements{We thank Matt Baker, Melody Chan, Alex Fink, Jeff Giansiracusa,  Felipe Rinc\'on, and Martin Ulirsch for helpful conversations. We particularly thank Chayim Lowen for pointing us to Eckmann-Hilton operators over hyperfields. Parts of this paper were written while V.S. was a Member at the Institute for Advanced Study, Princeton, and funded by the Charles Simonyi Endowment.}

\printbibliography

\clearpage
\section*{AI disclosure}

AI has not been used in the development of this manuscript at any stage.

\end{document}